\documentclass[12pt]{amsart}
\usepackage{amssymb}
\usepackage{amscd}

\newtheorem{Theorem}{Theorem}[section]
\newtheorem{Lemma}[Theorem]{Lemma}
\newtheorem{Corollary}[Theorem]{Corollary}
\newtheorem{Proposition}[Theorem]{Proposition}
\theoremstyle{definition}
\newtheorem{Definition}[Theorem]{Definition}
\newtheorem{Example}[Theorem]{Example}
\theoremstyle{remark}
\newtheorem{Remark}{Remark}

\font\sy=cmsy10

\font\ym=msbm10  
\newcommand{\cC}{{\hbox{\sy C}}}
\newcommand{\cD}{{\hbox{\sy D}}}

\newcommand{\cK}{{\hbox{\sy K}}}
\newcommand{\cL}{{\hbox{\sy L}}}

\newcommand{\cT}{{\hbox{\sy T}}}

\newcommand{\C}{{\text{\ym C}}}
\newcommand{\K}{{\text{\ym K}}}

\newcommand{\R}{\text{\ym R}}

\newcommand{\Z}{\text{\ym Z}}
\newcommand{\End}{\hbox{\rm End}}
\newcommand{\Hom}{\hbox{\rm Hom}}

\newcommand{\Ker}{\mathop{\rm Ker}}

\usepackage[dvips]{graphicx}
\newcommand{\trace}{\hbox{\rm tr}}

\title[]
{A Categorical and Diagrammatical 
Approach to Temperley-Lieb Algebras}
\author[Yamagami Shigeru]{Yamagami Shigeru}

\begin{document}
\maketitle   
\begin{center}
Department of Mathematical Sciences
\end{center}
\begin{center}
Ibaraki University 
\end{center}
\begin{center} 
Mito, 310-8512, JAPAN 
\end{center}    

\begin{abstract}
Algebraic basics on Temperley-Lieb algebras are proved 
in an elementary and straightforward way with the 
help of tensor categories behind them.
\end{abstract}
\bigskip
                           

\section{Introduction}
Temperley-Lieb algebra is a key notion in understanding 
quantum symmetry of various mathematical or physical objects 
and a variety of investigations have been worked out 
since the advent of the V.~Jones' 
celebrated work on knot invariants.
There have been developed three major approaches to the subject: 
the Jones' original method in subfactor theory 
together with the associated combinatorial invariants
(\cite{Jo}, \cite{GHJ}, \cite{EK}), 
representation theory of quantum groups 
(see \cite{CP}, \cite{Kas} for example), 
and the geometric (diagrammatic) method due to L.H.~Kauffman 
(\cite{Kau}, \cite{BJ}). 

In this paper, we shall present 
main structural analyses on Temperley-Lieb algebras, such as 
the criterion for semisimplicity, the description of Bratteli diagrams,  
the existence of C*-structures and so on, 
in a quite elementary and self-contained way 
with the emphasis on motivational streamlines of arguments. 
Our main tool is 
a diagrammatic presentation exploited in the third approach, 
together with naturally associated tensor categories 
of planar strings 
(called skein category in \cite{Tu} and 
Temperley-Lieb category in \cite{GW}), 
where Temperley-Lieb algebras are captured as 
the algebra of endomorphisms in Temperley-Lieb categories.

Widely recognized is the usefulness of 
such diagrammatic presentations in tensor calculus 
(see \cite{Tu}, for example) but planar strings themselves also
play substantial roles in describing quantum symmetries 
(see \cite{Kau2}, \cite{Tu} for example).

Returning to our main subject 
the Temperley-Lieb category belongs to the class of rigid tensor 
categories and if one manipulates its monoidal structure 
in a quite natural way, 
it already turns out to provide proofs of structural results mentioned 
above. 
The method is applicable to Fuss-Catalan algebras 
(\cite{BJ}, \cite{La}) as well 
(see \cite{FTC} for hints on the way of arguments) 
but we concentrate here just on 
the Temperley-Lieb case in viewing its fundamental importance 
among other related algebraic structures.  

Since the results themselves are well-known, we shall not repeat them 
here. 
Instead, we will briefly review existing approaches to the subject. 

Restricted to operator algebras, all the relevant analysis 
was worked out in \cite{Jo},
where 
the existence of (universal) Temperley-Lieb algebras 
is proved with the help of operator algebras of Murray and von Neumann. 
In that respect, the construction is highly analytical. 

On the other hand, representation theory of quantum groups 
has been mostly algebraic in its nature and, once one knows 
the relevant definitions, the whole analysis can be traced after 
the representation theory of ordinary compact Lie groups. 
One big conceptual gap here is the very definition of 
quantum groups, it was in fact a consequence
of ceaseless efforts of many researchers.

Compared to these, the approach here is straightforward, 
which is a combination of graphical presentation and 
rigidity calculus in tensor categories: 
the basic idea is just to try to determine the fusion rule 
(Clebsh-Gordan rule) in the Temperley-Lieb categories, 
so it is quite elementary up to topological intuition of 
planar isotopy of strings. 

The inductive formula for 
Jones-Wenzl idempotents are consequently derived as a byproduct of 
semisimplicity analysis. In this respect, our reasoning is reverse
in its order to standard arguments in 
\cite[Chap.~16]{Kau3}, \cite[Chap.~XII]{Tu}.  

The unitarity (positivity) criterion is also presented as a natural 
consequence of the present approach, which should be compared with 
the elaborate analysis in \cite{FK}. 



\section{Linear Categories}

By a {\bf linear category} we shall mean a category for which
hom-sets are vector spaces over a specified ground field $\K$
and all the relevant operations are assumed to be $\K$-linear. 
Therefore, given an object $X$ in a linear category, 
$\End(X) = \Hom(X,X)$ is a unital 
$\K$-algebra with the unit given by the identity morphism $1_X$. 
A linear category is said to be finite-dimensional if all hom-sets 
are finite-dimensional vector spaces. 

Since we have applications to physics or analysis in mind, 
the ground field is assumed to be the complex number field 
though the results can be formulated for a general field 
with possibly extra conditions depending on situations.


Recall that, in a linear category $\cL$, a direct sum 
$X_1\oplus \dots\oplus X_m$ 
is defined as an object $X$ together with morphisms 
$\alpha_j: X_j \to X$, $\beta_j:X \to X_j$ satisfying
$\beta_i\alpha_j = \delta_{i,j} 1_{X_i}$ and  
$1_X = \alpha_1\beta_1 + \dots + \alpha_n\beta_n$.
Each $X_i$ is called a direct summand of $X$. 
Since there could be many such morphisms, we shall often write
$X \cong X_1\oplus \dots\oplus X_m$ 
if we do not worry about their specific choices. 

If we take another direct sum $Y \cong Y_1\oplus \dots\oplus Y_n$, 
there arises the natural isomorphism of vector spaces 
\[
\Hom(X,Y) \cong \bigoplus_{i,j} \Hom(X_i,Y_j).
\]

Given a linear category $\cL$, its completion by idempotents
is, by definition, a linear category $\widetilde{\cL}$, 
where objects of $\widetilde{\cL}$ consist of a pair $(p,X)$
with $p \in \End(X)$ an idempotent and hom-sets are set to be
\[
\Hom((p,X),(q,Y)) = q\Hom(X,Y)p
\]
with the composition of morphisms given by the operation in $\cL$. 

We shall also use the notation $pX$ to stand for the object $(p,X)$ in 
$\widetilde\cL$. 

The notion of semisimplicity is usually formulated 
as a property on 
abelian categories. We shall deal with 
a specific class 
of non-abelian linear categories in what follows, for which we need to 
talk about the semisimplicity even though. 
We will here introduce it in a local and algebraic manner: 
First we extend a linear category $\cL$ by adding 
finite sequences of objects in $\cL$ with the notation 
$X_1\oplus \dots \oplus X_m$ to stand for the sequence 
$\{ X_1, \dots, X_m\}$ and the hom-sets among these are defined by 
\[
\Hom(X_1\oplus \dots \oplus X_m,Y_1 \oplus \dots \oplus Y_n) 
= 
\begin{pmatrix}
\Hom(X_1,Y_1) & \dots & \Hom(X_m,Y_1)\\
\vdots & \ddots & \vdots\\
\Hom(X_1,Y_n) & \dots & \Hom(X_m,Y_n)
\end{pmatrix}
\]
with the composition of morphisms given by matrix multiplication. 




\begin{Definition}
A finite family $\{ X_1, \dots, X_m\}$ of objects in 
a finite-dimensional linear category 
is said to be \textbf{semisimple} if the algebra 
\[
\End(X_1\oplus \dots \oplus X_m)
\]
is semisimple. 

A finite-dimensional linear category $\cL$ is said to be 
\textbf{essentially semisimle} 
if any finite family of objects in $\cL$ is semisimple. 
\end{Definition}

\begin{Lemma}
Let $A$ be a finite-dimensional semisimple algebra over the field 
$\C$. 
\begin{enumerate}
\item
Given an idempotent $p = p^2$ in $A$,  
we can find a finite family $\{ p_i\}$ of minimal idempotents 
such that 
\[
p = \sum_i p_i, 
\qquad
p_ip_j = \delta_{i,j} p_i, 
\]
which is referred to as a resolution of $p$ in $A$. 
\item 
For any resolution $\{ p_i\}$ of $p$, $p_iAp_j = \C \delta_{ij} p_i$. 
\end{enumerate}
\end{Lemma}

\begin{Corollary}
Given a finite-dimensional essentially semisimple $\C$-linear 
category $\cL$, we can find a family of objects $\{ X_i \}_{i \in I}$ 
in its idempotent-completion $\widetilde\cL$ such that 
$\Hom(X_i,X_j) = \C \delta_{ij} X_i$ and 
any object in $\widetilde\cL$ is isomorphic to a direct sum of 
finitely many objects in $\{ X_i\}$; 
\[
X \cong \bigoplus_{i \in I} m_iX_i
\quad
\text{with}
\ 
mX_i
= \overbrace{X_i\oplus \dots \oplus X_i}^{\text{$m$-times}}.
\]
Here the multiplicity function $m_i$ taking values in 
$\{ 0, 1, 2,\dots\}$ 
admits non-zero integers 
only on a finite subset of $I$. 

In other words, the linear category $\widetilde\cL$ is remade into 
a semisimple linear category by adding direct sums to $\widetilde\cL$. 
\end{Corollary}

\section{Tensor Categories}

By a tensor category, we shall here mean a finite-dimensional 
$\C$-linear monoidal category with the unit object $I$ 
satisfying $\End(I) = \C 1_I$. The strictness of monoidal structure 
is also assumed: $(f\otimes g)\otimes h = f\otimes (g\otimes h)$. 

Recall that, given a pair of objects $X$ and $Y$, 
$X$ is a left dual of $Y$ (or equivalently $Y$ is a right dual of $X$) 
if we can find a pair of morphisms $\epsilon: X\otimes Y \to I$ 
and $\delta: I \to Y\otimes X$ for which 
the following compositions are identities:
\[
\begin{CD}
X @>{1\otimes \delta}>> X\otimes Y\otimes X @>{\epsilon\otimes 1}>> X
\end{CD}
\]
\[
\begin{CD}
Y @>{\delta\otimes 1}>> Y\otimes X\otimes Y @>{1\otimes \epsilon}>> Y
\end{CD}
\]

The object $X$ (resp.~$Y$) is uniquely determined by $Y$ (resp.~$X$) 
up to isomorphisms and denoted by ${}^*Y$ (resp.~$X^*$). 
If we can choose $X = Y$, the object $X$ is said to be 
\textbf{self-dual}.  
A tensor category $\cT$ is said to be \textbf{rigid} if 
every object $X$ admits both left and right duals.

The operation of idempotent-completion is compatible with 
rigidity: 

\begin{Lemma}
Let $\cT$ be a rigid tensor category. Then its 
idempotent-completion 
$\widetilde{\cT}$ is rigid as well. 
\end{Lemma}

\begin{proof}
Let $Y$ be a right dual of an object $X$ with respect to 
morphisms $\epsilon: X\otimes Y \to I$ and $\delta: I \to Y\otimes X$. 
If $p \in \End(X)$ is an idempotent, 
the morphism $q \in \End(Y)$ defined by 
$q = (1_Y\otimes \epsilon)(1_Y\otimes p\otimes 1_Y)(\delta\otimes 1_Y)$
is an idempotent and $qY$ is a right dual of $pX$ by 
the morphisms 
\[
\epsilon(p\otimes q) = \epsilon(p\otimes 1_Y) 
: pX\otimes qY \to I,
\quad 
(q\otimes p)\delta = (1_Y\otimes p)\delta: I \to qY\otimes pX.
\]
\end{proof}

The next is an immediate consequence of rigidity as is well-known. 

\begin{Lemma}[Frobenius Reciprocity]
In a tensor category $\cT$, if an object $X$ admits a right dual $X^*$, 
then we have the natural isomorphism 
\[
\Hom(Y,Z\otimes X) \cong \Hom(Y\otimes X^*,Z).
\]
\end{Lemma}

Let $X$ and $Y$ be objects in a tensor category. 
Then the map $\End(X) \ni a \mapsto a\otimes 1_Y \in \End(X\otimes Y)$ 
is a unital homomorphism. To see the injectivity of this map, 
we observe the following: 

\begin{Lemma}
If $Y$ admits a right dual $Y^*$ such that 
$Y\otimes Y^* \cong I \oplus Z$ for some object $Z$, i.e., 
we can find morphism $\delta': I \to Y\otimes Y^*$ satisfying 
$\epsilon \delta' = 1_I$, 
then the algebra-homomorphism $\End(X) \to \End(X\otimes Y)$ is 
injective. 
\end{Lemma}

Recall that, given an inclusion $A \subset B$ of finite-dimensional 
semisimple algebras with the common unit element, 
its Bratteli diagram is a bipartite graph whose vertex set is 
the disjoint union of the set $\widehat A$ of 
equivalence classes of simple $A$-modules and 
the set $\widehat B$ of equivalence classes of simple $B$-modules 
with two vertices $i \widehat A$ and $j \in \widehat B$ 
connected by $m$-edges, 
where a non-negative integer $m$ 
is determined as follows: letting ${}_AX$ and ${}_BY$ be simple modules
representing $i$ and $j$, we set 
\[
m = \dim \Hom({}_AX, {}_AY).
\]

For the inclusion of algebras in the previous lemma, the following 
is, though immediate, fundamental. 
 

\begin{Lemma}
Under the same assumption as in the above lemma, assume that
there are objects $\{ X_i\}$ and $\{ Y_k\}$ satisfying
\begin{gather*}
\Hom(X_i,X_j) = \delta_{i,j}\C 1_{X_i}, 
\qquad 
\Hom(Y_k,Y_l) = \delta_{k,l}\C 1_{Y_k},\\
\Hom(X_i,Y_k) = \{ 0\} = \Hom(Y_k,X_i)
\end{gather*}
and 
\[
X_i\otimes Y \cong \bigoplus_j m_{ij} X_j
\oplus 
\bigoplus_k n_{ik} Y_k. 
\]
Then, for an object $X$ which is isomorphic to 
a direct sum of finitely many $X_i$'s, both of $\End(X)$ and 
$\End(X\otimes Y)$ are semisimple with the Bratteli diagram 
of the inclusion $\End(X) \subset \End(X\otimes Y)$ 
specified by $m_{ij}$ and $n_{ik}$ as the numbers of edges. 
\end{Lemma}

\section{Temperley-Lieb Categories}


Here we shall review the planar description of 
Temperley-Lieb algebras
according to L.H.~Kauffman (\cite{Kau}, \cite{Kau2}), 
which leads us to the accompanied tensor categories 
at the same time (cf.~\cite{Tu}, \cite{Ba}, \cite{GW}). 



Let $m$ and $n$ be non-negative integers of the same parity. 
Choose a rectangle in the plane with marking of $m$ points on 
the upper horizontal edge 
and $n$ points on the lower horizontal edge. 
Join these $m+n$ points by $(m+n)/2$ planar curves 
inside the rectangular box 
so that curves do not cross each other. 
We denote by $K_{m,n}$ the set of
isotopy classes of planar curves of this type. 
An element in $K_{m,n}$ is referred to as a Kauffman diagram of 
type $(m,n)$. 
The set $K_{n,n}$ is simply denoted by $K_n$. 
It is well-known that the set $K_{m,n}$ consists of 
$C_{(m+n)/2}$ diagrams, where $C_n = \binom{2n}{n}/(n+1)$ 
denotes the $n$-th Catalan number. 
When $m$ and $n$ have different parity, 
we set $K_{m,n} = \emptyset$. 


See Figure \ref{three} 
(the bounding boxes being omitted) 
for the isotopy patterns in $K_3 = K_{3,3}$.

\bigskip
\begin{figure}[h]
\input{three.tpc}
\caption{\label{three}}
\end{figure}
\bigskip

Let $\C[K_{m,n}]$ be the free $\C$-vector space generated by 
the set $K_{m,n}$, which is also denoted by $\cK_d(m,n)$. 
Given a complex number $d$, 
a linear category $\cK_d$ is defined in the following way: 
objects of $\cK_d$ are non-negative integers and hom-sets are 
set to be $\Hom(n,m) = \C[K_{m,n}] = \cK_d(m,n)$ with the composition 
of morphisms given by 
the concatenation of planar strings 
through one horizontal edge of boxes 
(diagrams stream from bottom to top by convention), 
where each loop (if there appeared any) is 
(removed and) replaced by the complex 
number $d$ (Figure \ref{loop}).


\bigskip
\begin{figure}[h]
\hspace{1cm}
\input{loop.tpc}
\caption{\label{loop}}
\end{figure}

Taking the dependence on $d$ into account, we also use the notation 
$\C[K_{m,n},d]$. 

For $D \in K_{m,n}$ and $D' \in K_{m',n'}$, 
we define $D\otimes D' \in K_{m+m',n+n'}$ 
by placing $D$ and $D'$ horizontally 
so that $D$ is left to $D'$ 
(juxtaposition). 
The operation is clearly associative and is linearly extended to 
the map $\C[K_{m,n}]\otimes \C[K_{m',n'}] \to \C[K_{m+m'n+n'}]$, 
which makes $\cK_d$ into a tensor category, 
called the \textbf{Temperley-Lieb category}. 
The terminology is in accordance with \cite{Ba}, \cite{GW} 
though it is called skein category in \cite{Tu}.)

In the category $\cK_d$, 
the multiplicative notation is also used to indicate objects; 
we introduce a dummy symbol $X$ to represent the object $1$ so that 
the object $n$ is expressed by 
\[
X^{\otimes n} = 
\overbrace{X\otimes \dots \otimes X}^{\text{$n$-times}}.
\]
For short, we use the notation $X^n$ occasionally.

The Temperley-Lieb category possesses 
a specific feature of perfect rigidity: 
Firstly, the generating object $X$ is self-dual with respect to 
the pairing $\epsilon: X\otimes X \to I$ and 
the copairing $\delta:I \to X\otimes X$ given by the arcs 
(Figure \ref{arcs}). 

\begin{figure}[h]
\hspace{1cm}
\input{arcso.tpc}
\caption{\label{arcs}}
\end{figure}

Iterating these basic morphisms, we see that every object 
is self-dual: The pairing and coparing of $X^3$, for example, 
are given by the diagrams in Figure \ref{multarcs}. 

\begin{figure}[h]
\hspace{1cm}
\input{multarcs.tpc}
\caption{\label{multarcs}}
\end{figure}

Frobenius transforms are then visually realized 
by bending terminal lines so that it changes directions of morphisms
(Figure \ref{bending}). 

\begin{figure}[h]
\input{bending.tpc}
\caption{\label{bending}}
\end{figure}

As an application of this geometrical interpretation, 
we see that two ways of complete bending coincide
(Figure \ref{transposed}).

\begin{figure}[h]
\hspace{-5mm}
\input{transposed.tpc}
\caption{\label{transposed}}
\end{figure}

In fact, 
given a diagram $D \in K_{m,n}$, both of these operations
result in the transposed diagram ${}^tD \in K_{n,m}$ 
which is, by definition, the rotation of $D$ by 
an angle of $\pi$. 
See Figure~\ref{rotation} as an example of $D \in K_{3,1}$ and 
${}^tD \in K_{1,3}$. 

\begin{figure}[h]
\hspace{2cm}
\input{rotation.tpc}
\caption{\label{rotation}}
\end{figure}

The operation is linearly extended to the map 
$\Hom(X^m,X^n) \to \Hom(X^n,X^m)$, which is clearly involutive and 
antimultiplicative: ${}^t({}^tf) = f$ and 
${}^t(fg) = ({}^tg)({}^tf)$. Moreover, it is compatible 
with the monoidal structure in the sense that 
${}^t(f\otimes g) = {}^tg\otimes {}^tf$ for 
$f \in \Hom(X^m,X^{m'})$ and $g \in \Hom(X^n,X^{n'})$, 
i.e., the duality holds for rigidity. 

On each algebra $\End(X^n)$, a special functional $\trace_n$ 
(called \textbf{Markov trace})
is defined by closing diagrams completely on $K_n$ and then 
taking the linear extension to $\End(X^n) = \C[K_n]$. 
Apparently we have two choices 
for closing in the plane, which, however, gives the same result 
as indicated by Figure \ref{trace}. 
The suffix $n$ is often omitted if it causes no confusion. 
(See \cite{Mal}, \cite{Da}, \cite{BW} for generalities  
on duality and traces.)


\begin{figure}[h]
\input{trace.tpc}
\caption{\label{trace}}
\end{figure}

Here are some of formulas concerning the Markov trace
and transposed morphisms: 
\begin{enumerate}
\item
$\trace_n(fg) = \trace_m(gf)$ for 
$f \in \Hom(X^m,X^n)$ and $g \in \Hom(X^n,X^m)$.
\item
$\trace_n(h) = \trace_n({}^th)$ for $h \in \End(X^n)$.
\end{enumerate}


\begin{Remark}
The reflection of diagrams vertically or horizontally gives 
another involution on hom-sets. These with the rotation as  
transposed morphisms constitute a symmetry of $\Z_2\times \Z_2$. 
\end{Remark}

Let us now introduce \textbf{elementary diagrams} 
$h_1, \dots, h_{n-1}$ in $K_n$ by Figure \ref{jones}

\begin{figure}[h]
\hspace{-1cm}
\input{jones.tpc}
\caption{\label{jones}}
\end{figure}

\noindent 
which satisfy the relations of Temperley-Lieb 
in $\End(X^n) = \C[K_n,d]$:
\[
h_i^2 = d h_i, 
\quad 
h_ih_j = h_jh_i 
\quad (|i-j| \geq 2), 
\quad 
h_i h_{i\pm 1}h_i = h_i.
\]

By adding one vertical line (with two end points) to 
the right end of planar strings in $K_n$, we have the imbedding 
$K_n \subset K_{n+1}$, which induces an inclusion of algebras 
$\C[K_n,d] \subset \C[K_{n+1},d]$. 
In terms of the monoidal structure, this is expressed by 
$a \mapsto a\otimes 1$.

The elementary diagrams $\{ h_1, \dots, h_{n-1} \}$, 
together with the unit diagram $1$, 
turn out to generate the algebra $\C[K_n,d]$ 
(cf.~\cite[Theorem 4.3]{Kau2}):
In fact, given a diagram $D$ in $K_n$, we stretch out each 
string in $D$ vertically and then wave it horizontally. 
Then many minimal arcs are coupled by cutting strings horizontally at 
levels without critical crossings, resulting in a product formula 
for a planar diagram in terms of the elementary diagrams 
$\{ h_1, \dots, h_{n-1} \}$ (Figure \ref{morse}). 
(See Appendix~A for a more detailed account.)

\begin{figure}[h]
\input{morseo.tpc}
\caption{\label{morse}}
\end{figure}

The following can be read off from Figure \ref{iterate}, 
which already gives a pictorial proof of the formula 
for iterated basic constructions discussed in \cite{PP} 
(see \cite{NOA} for the tensor-categorical meaning of 
the Jones basic construction).

\begin{Example}
The $h$-element for the object $X^{\otimes n}$ is given by 
\[
(h_nh_{n-1}\dots h_1)(h_{n+1}h_n\dots h_2)\dots 
(h_{2n-1}h_{2n-2}\dots h_n).
\]
\end{Example}

\begin{figure}[h]
\input{iterateo.tpc}
\caption{\label{iterate}}
\end{figure}

The above observation also reveals the fact that hom-sets are generated 
by basic arcs $\epsilon:X\otimes X \to I$ and 
$\delta: X\otimes X \to I$ together with their tensor product 
ampliations $1_{X^m}\otimes \epsilon\otimes 1_{X^n}$ and 
$1_{X^m}\otimes \delta\otimes 1_{X^n}$. 

In particular, a monoidal functor $F$ on the tensor category 
$\cK_d$ is uniquely determined by morphisms $F(\epsilon)$ and 
$F(\delta)$. 
When $F: \cK_d \to \cK_{d'}$, we should have 
$F(\epsilon) = \lambda \epsilon$ and 
$F(\delta) = \mu \delta$ with $\lambda, \mu \in \C^\times$. 
Since $F$ must preserve hook identities, we are forced to 
set $\lambda\mu = 1$, which in turn ensures the whole multiplicativity 
of $F$. In particular, we should have 
\[
d'1_I = F(\epsilon)F(\delta) 
= F(\epsilon\delta) = F(d 1_I) = d1_I.
\]
Thus Temperle-Lieb categories $\cK_d$ for different $d$ 
are not equivalent as tensor categories. 

The commutation relations of Temperley-Lieb algebra originally 
emerged in a model of statistical physics (\cite{TL}), 
which were later rediscovered by 
V.~Jones as commutation relations among idempotents (\cite{Jo}). 

If we set $e_i = h_i/d$, these are idempotents and the Temperley-Lieb 
relations are equivalently described by the relations 
\[
\begin{cases}
  e_ie_je_i = d^{-2} e_i &\text{if $|i-j| = 1$,}\\ 
  e_ie_j = e_je_i &\text{if $|i-j| \geq 2$.}
\end{cases} 
\]
We remark here that this simple observation shows that 
$\C[K_n,d_1] \cong \C[K_n,d_2]$ if $d_1^2 = d_2^2$. 

Given an integer $n \geq 1$, let $A_n$ be the algebra 
universally generated by 
$\{ h_1, \dots, h_{n-1} \}$ 
and the unit $1$ 
with the Temperley-Lieb relations, which is referred to as 
the \textbf{Temperley-Lieb algebra}. 
(Rigorously speaking, we should use other symbols, 
say $h_i'$, instead of $h_i$.)
By universality, we have the natural epimorphism $A_n \to \C[K_n]$, 
which turns out to be an isomorphism. 
(We reproduce somewhat simplified proofs in Appendix A.)  

\begin{Proposition}[{\cite[Theorem 4.3]{Kau2}}]
The algebra $\C[K_n,d]$ is 
universally generated by 
$\{ h_1, \dots, h_{n-1} \}$ and the unit $1$ with the relations 
\[
h_i^2 = d h_i, 
\quad 
h_ih_j = h_jh_i 
\quad (|i-j| \geq 2), 
\quad 
h_i h_{i\pm 1}h_i = h_i
\]
for $1 \leq i,j \leq n-1$, whence it is identified with 
the Temperley-Lieb algebra $A_n$. 

In particular, for $n \geq 1$, 
the obvious homomorphism $A_n \to A_{n+1}$ of 
Temperley-Lieb algebras is injective.   
\end{Proposition}

\bigskip

\section{Semisimplicity and Fusion Rules}

Here we shall analyse the semisimplicity of the Temperley-Lieb 
categories together with their fusion rules, by looking into 
the structure of the algebra $A_n = \End(X^n)$ for 
$n \geq 1$. 

Since $A_2 = \langle 1, h_1\rangle$ and $h_1^2 = d h_1$, 
the algebra $A_2$ is semisimple if and only if $d \not= 0$ and, 
if this is the case, 
$A_2 = \C e_1 \oplus \C(1-e_1) \cong 
\C \oplus \C$ with $e_1 = h_1/d$ an idempotent. In other words, 
the linear subcategory generated by $\{ I, X, X^2\}$ is semisimple, 
with simple objects given by $I$, $X$ and $X_2 \equiv f_2X^2$, where 
$f_2 = 1-e_1$ and $e_1X^2 \cong I$; $X^2 \cong I \oplus X_2$. 
Note also that $X_2$ is self-dual because ${}^tf_2 = f_2$. 
For the notational consistency, we also write $f_1 = 1_X$ and 
$X_1 = f_1X = X$. 

Next, under the assumption of semisimplicity at the first stage, 
i.e., $d \not= 0$, we see 
$X^3 \cong (I\oplus X_2)\otimes X \cong X \oplus X_2\otimes X$ 
in the tensor category $\widetilde{\cK}_d$. 
So we need to investigate how $X_2\otimes X$ is interrelated to $X$. 
By Frobenius reciprocity, we have 
\begin{align*}
\Hom(X,X_2\otimes X) &\cong \Hom(X\otimes X,X_2) 
\cong \Hom(I\oplus X_2,X_2) = \End(X_2),\\
\Hom(X_2\otimes X,X) &\cong \Hom(X_2,X\otimes X) 
\cong \Hom(X_2,I\oplus X_2) = \End(X_2).
\end{align*}

Let $\varphi_2: X \to X_2\otimes X$ and 
$\psi_2: X_2\otimes X \to X$ be Frobenius transforms of 
$1_{X_2} = f_2$ through the above isomorphisms 
(see Figure \ref{phipsi}). 
Then any morphism $X \to X_2\otimes X$ is a scalar multiple of 
$\varphi_2$ and similarly for $\psi_2$. 

By a diagrammatic computation (Figure \ref{critical}), we see that 
$\psi_2\varphi_2 = \lambda_2 1_X$ with $\lambda_2 = d - d^{-1}$. 
Thus, if $\lambda_2 \not= 0$, 
$X$ is a direct summand in $X_2\otimes X$ 
with the projection to $X$ given by
$\lambda_2^{-1} \varphi_2\psi_2 \in \End(X_2\otimes X)$ and 
the complementary subobject $X_3$ of $X_2\otimes X$ given by 
the idempotent $f_3 = f_2\otimes 1_X - \lambda_2^{-1}\varphi_2\psi_2$. 
Since $f_3$ is an idempotent in 
$\End(X_2\otimes X) = (f_2\otimes 1_X)\End(X^3)(f_2\otimes 1_X)$, 
it is also an idempotent in $A_3 = \End(X^3)$ 
(satisfying $f_3(f_2\otimes 1_X) = (f_2\otimes 1_X)f_3 = f_3$.).

Since both of $\Hom(X,X_2\otimes X) \cong \Hom(X, X \oplus X_3)$ and 
$\Hom(X_2\otimes X,X) \cong \Hom(X\oplus X_3,X)$ are one-dimensional, 
we see that 
\[
\Hom(X,X_3) = \Hom(X_3,X) = \{ 0\},
\]
while the parity condition implies the triviality of 
$\Hom(I,X_3)$, $\Hom(X_3,I)$, $\Hom(X_2,X_3)$ and $\Hom(X_3,X_2)$. 
Note here that 
$X^3 \cong 2X \oplus X_3$. 

If $f_3 = 0$, i.e., $X_3 = 0$, then 
$\End(X^3) \cong M_2(\C)$ is four-dimensional, a contradiction 
with $\dim \End(X^3) = 5$. 
Thus $f_3 \not= 0$ and 
the dimension estimate 
\[
2^2 + 1 \leq \dim \End(X\oplus X) + \dim \End(X_3) = \dim A_3 = 5
\]
shows $\End(X_3) = \C f_3$. 


\begin{figure}[h]
\input{critical.tpc}
\caption{\label{critical}}
\end{figure}

Summarizing the discussion so far, 
under the assumption $d \not= 0$ and $d - d^{-1} \not= 0$, 
the linear subcategory of $\cK_d$ generated by 
$\{ I, X, X^2, X^3\}$ is semisimple with the (isomorphism classes of) 
simple objects given by $I$, $X$, $X_2$ and $X_3$. 

The reasoning is applicable repeatedly and we arrive at 
the following induction scheme: 
Assume that idempotents $f_k \in A_k$ are inductively 
defined up to $k = n$ so that 
\begin{enumerate}
\item the sequence $\{ f_k\}$ satisfies the recursive formula 
\[
f_{k+1} = f_k\otimes 1_X - 
\frac{\trace(f_{k-1})}{\trace(f_k)} 
(f_k\otimes 1_X) h_k (f_k\otimes 1_X)
\]
with $\trace(f_k) \not= 0$ for $1 \leq k < n$
\item 
the linear subcategory generated by $\{ I, X, X^2, \dots, X^n\}$ 
is semisimple with inequivalent simple objects 
represented by $X_k \equiv f_kX^{\otimes k}$ for $0 \leq k \leq n$, 
\item 
$X_k\otimes X \cong X_{k-1} \oplus X_{k+1}$ for $1 \leq k <n$ 
with $X_0 = I$. 
\end{enumerate}

At this stage, we derive two consequences from the above hypotheses: 
Applying the Markov trace to the recursive formula for $f_j$, 
we have 
\[
\trace(f_{k+1}) = 
d\, \trace(f_k) - \trace(f_{k-1})
\]
for $1 \leq k <n$ with $\trace(f_0) =1$ and 
$\trace(f_1) = d$. 
Consequently, by the choice $d = q + q^{-1}$ with $0 \not= q \in \C$, 
we have 
\[
\trace(f_k) = [k+1]_q =  
\frac{q^{k+1} - q^{-k-1}}{q - q^{-1}}.
\]

From the fusion rule for $X_k\otimes X$, 
together with lattice path countings 
by the reflection technique, 
we obtain the multiplicity formula (see Figure \ref{ladder})
\[
X^{\otimes k} = \bigoplus_{j=0}^{[k/2]} 
\begin{bmatrix}
k\\
j
\end{bmatrix}
X_{k-j}
\]
for $1 \leq k \leq n$, where 
\[
\begin{bmatrix}
k\\
j
\end{bmatrix} 
= \binom{k}{j} - \binom{k}{j-1}.
\]
($\binom{k}{-1} = 0$ by definition.)

Here are explicit formula in lower cases: 
\begin{align*}
X^2 &= X_2 \oplus I,\\
X^3 &= X_3 \oplus 2X,\\
X^4 &= X_4 \oplus 3X_2 \oplus 2I,\\
X^5 &= X_5 \oplus 4X_3 \oplus 5X.
\end{align*}

\begin{figure}
\input{ladder.tpc}
\caption{\label{ladder}}
\end{figure}

We can now raise the induction stage one step further: 
By Frobenius reciprocity and the induction hypothesis on 
the fusion rule,  we have
\begin{align*}
\Hom(X_k,X_n\otimes X) 
&\cong \Hom(X_{k-1} \oplus X_{k+1},X_n) = \{0\},\\
\Hom(X_n\otimes X,X_k)
&\cong \Hom(X_n,X_{k-1} \oplus X_{k+1}) = \{ 0\}
\end{align*}
for $1 \leq k \leq n-2$, which means that no $X_k$-component 
appears in $X_n\otimes X$ for $1 \leq k \leq n-2$. 
On the other hand, again by Frobenius reciprocity and the fusion rule 
assumption, 
\begin{align*}
\Hom(X_{n-1},X_n\otimes X) &\cong \Hom(X_{n-1}\otimes X,X_n) 
\cong \Hom(X_{n-2}\oplus X_n,X_n) = \End(X_n),\\
\Hom(X_n\otimes X,X_{n-1}) &\cong \Hom(X_n,X_{n-1}\otimes X) 
\cong \Hom(X_n,X_{n-2}\oplus X_n) = \End(X_n).
\end{align*}
Thus, as Frobenius transforms of $1_{X_n} = f_n$, we can define 
non-zero morphisms 
$\varphi_n: X_{n-1} \to X_n\otimes X$ and 
$\psi_n: X_n\otimes X \to X_{n-1}$ (Figure \ref{phipsi}).

\begin{figure}[h]
\input{psiphi.tpc}
\caption{\label{phipsi}}
\end{figure}

\noindent
Then, by manipulating diagrams 
(Figure \ref{lambda}, Figure \ref{markov}), we see 
$\psi_n\varphi_n = \lambda_n 1_{X_n}$ 
with 
$\lambda_n = \trace(f_n)/\trace(f_{n-1})$. 

\begin{figure}[h]
\input{lambda.tpc}   
  \caption{\label{lambda}}
\end{figure}

\begin{figure}[h]
\input{markov.tpc}
\caption{\label{markov}}
\end{figure}

\noindent
Therefore, if $\trace(f_n) \not= 0$, then we can define an 
idempotent $f_{n+1} \in A_{n+1}$ by the formula
\[
f_{n+1} = 
f_n\otimes 1_X - 
\frac{\trace(f_{n-1})}{\trace(f_n)} \varphi_n\psi_n 
= f_n\otimes 1_X - 
\frac{\trace(f_{n-1})}{\trace(f_n)}
(f_n\otimes 1_X)h_n(f_n\otimes 1_X)
\]
(Figure \ref{oldstuff}) with the associated subobject 
$X_{n+1} = f_{n+1} X^{\otimes(n+1)}$ and we reach the direct sum 
decomposition 
$X_n\otimes X \cong X_{n-1} \oplus X_{n+1}$. 
Since both of $\Hom(X_{n-1},X_n\otimes X)$ and 
$\Hom(X_n\otimes X,X_{n-1})$ are one-dimensional, we have 
\[
\Hom(X_{n-1},X_{n+1}) = \{ 0\} = \Hom(X_{n+1},X_{n-1}),
\]
whereas the triviality of $\Hom(X_n,X_{n+1})$ and 
$\Hom(X_{n+1},X_n\otimes X)$ is a consequence of parity discrepancy. 

\begin{figure}[h]
\input{oldstuff.tpc}
\caption{\label{oldstuff}}
\end{figure}

From the multiplicity formula for $X^{\otimes n}$ and 
the fusion rule $X_n\otimes X \cong X_{n-1} \oplus X_{n+1}$ 
with $\Hom(X_k,X_{n+1}) = \{ 0\} = \Hom(X_{n+1},X_k)$ for 
$1 \leq k \leq n$, we obtain the decomposition 
\[
X^{\otimes (n+1)} = \bigoplus_{j=0}^{[(n+1)/2]} 
\begin{bmatrix}
n+1\\
j
\end{bmatrix}
X_{n+1-j}.
\]

At this point, we have no information on
the simplicity of the new stuff $X_{n+1}$ yet. 
The above decomposition, however, gives rise to 
the following dimension identity
\[
\dim \End(X^{\otimes(n+1)}) = 
\dim\End(X_{n+1}) - 1 +  
\sum_{j=0}^{[(n+1)/2]} 
\begin{bmatrix}
n+1\\
j
\end{bmatrix}^2,
\]
which particularly implies $f_{n+1} \not= 0$. 

Now we conclude $\End(X_{n+1}) = \C f_{n+1}$ 
from the combinatorial formula below, which is obtained 
by folding halfway in the following well-known binomial identity
(see \cite[Chapter 5]{GKP} for example)
\[
\sum_k \binom{m}{a+k} \binom{n}{b+k} 
= \binom{m+n}{m-a+b} = 
\binom{m+n}{n+a-b}.
\]

\begin{Lemma}[\cite{Jo, GHJ}]
For a positive integer $n$, we have 
\[
\sum_{j=0}^{[n/2]}\, 
\begin{bmatrix}
n\\
j
\end{bmatrix}^2 
= \frac{1}{n+1} \binom{2n}{n}. 
\]
\end{Lemma}

For $n=5$, this means
\[
1^2 + 4^2 + 5^2 = 42.
\]

As a conclusion of induction arguments so far, we have 

\begin{Proposition}
Express $d$ in the form $d = q + q^{-1}$ with 
$q \in \C^\times$. Assume that $[k]_q \not= 0$ for 
$1 \leq k \leq n$. Then the linear subcategory of $\cK_d$ 
generated by 
$\{ I, X, \dots, X^{\otimes n} \}$ is semisimple and 
a representative simple objects $X_k$ ($k=1,2, \dots, n$) 
in $\widetilde{\cK}_d$ 
are inductively defined so that 
$X_k\otimes X \cong X_{k-1} \oplus X_{k+1}$ ($1 \leq k <n$). 

The subobject $X_k$ appears only once in $X^{\otimes k}$ and 
the associated idempotent $f_k \in \End(X^{\otimes k})$ is 
inductively defined by the Wenzl's formula
\[
f_{k+1} = f_k\otimes 1_X - 
\frac{[k]_q}{[k+1]_q} 
(f_k\otimes 1_X) h_k (f_k\otimes 1_X)
\]
with the trace value given by $\trace(f_k) = [k+1]_q$. 
\end{Proposition}

\begin{Corollary}
If $q^2$ is not a proper root of unity, the Temperley-Lieb category 
$\cK_d$ is essentially semisimple with its fusion rule given by 
the Clebsh-Gordan rule 
\[
X_j\otimes X_k \cong X_{|j-k|} \oplus X_{|j-k|+2} 
\oplus \dots \oplus X_{j+k}.
\]
\end{Corollary}


If $q^2$ is an $l$-th primitive root of unity, the non-degeneracy 
of the Markov trace 
ceases at the algebra $\End(X^{\otimes (l-1)})$. 
It is then customary to take a quotient of the tensor category 
$\cK_d$: Let 
\[
\Ker(X^{\otimes m},X^{\otimes n}) = 
\{ f \in \Hom(X^{\otimes m}, X^{\otimes n}); 
\langle fg \rangle = 0\  
\text{for $g \in \Hom(X^{\otimes n}, X^{\otimes m})$}
\}.
\]
Then these constitute an ideal of $\cK_d$ in the sense that 
\begin{enumerate}
\item
$\Hom(X^{\otimes m},X^{\otimes n}) 
\Ker(X^{\otimes l},X^{\otimes m})
\Hom(X^{\otimes k},X^{\otimes l}) \subset 
\Ker(X^{\otimes k},X^{\otimes n})$ and 
\item 
$\Hom(X^{\otimes m'},X^{\otimes n'})
\otimes \Ker(X^{\otimes m},X^{\otimes n}) 
\otimes \Hom(X^{\otimes m''},X^{\otimes n''}) 
\subset \Ker(X^{\otimes(m'+m+m'')},X^{\otimes(n'+n+n'')})$, 
\end{enumerate}
where the former is a consequence of the trace property of 
$\trace(\cdot)$ and the latter is checked by means of 
the conditional expectation $\End(X^{m+n}) \to \End(X^m)$ 
with respect to the inclusion 
$\End(X^m) \cong \End(X^m)\otimes 1 \subset \End(X^{m+n})$. 

The quotient tensor category $\overline{\cK}_d$ is then defined 
so that objects are the same with those for $\cK_d$ 
(but we shall use the bar notation to indicate objects in 
$\overline{\cK}_d$) and hom-sets are given by 
$\Hom(\overline{Y},\overline{Z}) = \Hom(Y,Z)/\Ker(Y,Z)$ 
with
the monoidal structure on $\overline{\cK}_d$ 
inherited from $\cK_d$.
The resultant tensor category is then referred to as 
the \textbf{reduced Temperley-Lieb category}. 

Now the following is immediate from our discussions so far. 

\begin{Proposition}
Assume that $q^2$ is an $l$-th primitive root of unity.
Then 
\[
\Ker(X^{\otimes k}) = 
\begin{cases}
\{ 0\} & \text{if $k < l-1$,}\\
\C f_{l-1} & \text{if $k = l-1$}
\end{cases}
\]
and $\Ker(X^{\otimes m},X^{\otimes n})$ is 
monoidally generated by $f_{l-1}$. 

The reduced Temperley-Lieb category is semisimple 
with simple objects given by 
$\{ \overline{X}_k\}_{0 \leq k \leq l-2}$ 
with the recursive formula
\[
\overline{X}_k\otimes \overline{X} = 
\begin{cases}
\overline{X}_{k-1}\oplus \overline{X}_{k+1} &\text{if $k < l-2$,}\\
\overline{X}_{l-3} &\text{if $k= l-2$.}
\end{cases}
\]
\end{Proposition}

\begin{Corollary}
The fusion rule for $\{ \overline{X}_k \}$ is given by 
the truncated Clebsh-Gordan rule:
\[
\overline{X}_j\otimes \overline{X}_k 
\cong \overline{X}_{|j-k|} \oplus \overline{X}_{|j-k|+2} \oplus \dots 
\oplus \overline{X}_m,
\]
where 
\[
m = 
\begin{cases}
j+k &\text{if $j+k \leq l-2$,}\\
2(l-2) -(j+k) &\text{if $j+k \geq l-2$.}
\end{cases}
\]
\end{Corollary}

\begin{Remark}
The kernel of $\cK_d$ is characterized in \cite{GW} 
as the unique monoidal ideal.
\end{Remark}

\section{Positivity in Temperley-Lieb Categories}

We shall now clarify the condition when the Temperly-Lieb 
category is a C*-tensor category. 
Recall that a linear category $\cL$ is a \textbf{C*-category} if 
hom-sets are Banach spaces and we are given conjugate linear 
maps (denoted by * and referred to as a star operation) 
on hom-sets 
$\Hom(Y,Z) \ni f \mapsto f^* \in \Hom(Z,Y)$ 
satisfying 
(i) $(f^*)^* = f$, (ii) $(fg)^* = g^*f^*$, and 
(iii) $\| f^*f\| = \| f\|^2$ for $g: X \to Y$, $f:Y \to Z$. 
When hom-sets are finite-dimensional, a more algebraic formulation 
is possible (see \cite[Appendix]{GHJ}): given a star operation 
satisfying (i), (ii) and the condition that $f^*f = 0$ implies $f=0$,  
there is the unique C*-norm fulfilling (iii). 

A \textbf{C*-tensor category} (or tensor C*-category) is 
a (strict) tensor category 
which is a C*-category at the same time with the common underlying 
linear structure such that two structures are compatible 
in the sense that 
$(f\otimes g)^* = f^*\otimes g^*$ for morphisms $f$, $g$. 
(when the associativity transformations are explicit, they are 
assumed to be unitary with respect to the star operation). 

A functor $F: \cC \to \cD$ between two (strict) tensor categories 
is said to be \textbf{monoidal} if $F(I_\cC) = I_{\cD}$,  
$F(X\otimes Y) = F(X)\otimes F(Y)$ for objects $X$, $Y$ and 
$F(f\otimes g) = F(f)\otimes F(g)$ for morphisms $f$, $g$.

Two tensor categories are said to be 
\textbf{monoidally equivalent} 
if we can find a monoidal functor between these tensor categories 
which gives an equivalence of categories; 
if the functor $F$ is 
fully faithful in the sense that $F$ gives isomorphisms 
on hom-sets and any object of $\cD$ is isomorphic to 
$F(X)$ for some object $X$ in $\cC$. 

When $\cC$ and $\cD$ are C*-tensor categories, a monoidal 
functor is said to be \textbf{C*-monoidal} if $F(f)^* = F(f^*)$ 
for any morphism $f$ in $\cC$. 

\begin{Remark}
Our definition of monoidality is the one usually referred to 
as being strict. Since any non-strict monoidal functor is 
changed to be strict by replacing tensor categories 
with equivalent ones, there are no essential differences. 
\end{Remark}

The following reveals the universality of Temperley-Lieb categories 
concerning self-dual objects in tensor categories.

\begin{Lemma}
Let 
$Y$ be a self-dual object in 
a tensor category $\cT$ 
with the associated morphisms $\epsilon_Y: Y\otimes Y\to I$, 
$\delta_Y:I \to Y\otimes Y$ and suppose that 
$\epsilon_Y\delta_Y = d 1_I$ with $d \in \C^\times$. 

Then the correspondence $\epsilon \mapsto \epsilon_Y$, 
$\delta \mapsto \delta_Y$ is extended to a monoidal functor 
from the Temperley-Lieb category $\cK_d$ to $\cT$. 
\end{Lemma}

\begin{proof}
By a single arc, we shall mean a morphism of the form 
$1\otimes \epsilon\otimes 1$ or $1\otimes \delta\otimes 1$. 
Given an arc $a$ in $\cK_d$, we denote by $F(a)$ the morphism 
in $\cT$ defined by 
\[
F(a) = 
\begin{cases}
1_{Y^m}\otimes \epsilon_Y\otimes 1_{Y^n} 
&\text{if $a = 1_{X^m}\otimes \epsilon \otimes 1_{X^n}$,}\\
1_{Y^m}\otimes \delta_Y\otimes 1_{Y^n} 
&\text{if $a = 1_{X^m}\otimes \delta \otimes 1_{X^n}$.}
\end{cases}
\]

As in the proof of the generating property of elementary diagrams, 
we see that each diagram $D \in K_{m,n}$ can be expressed 
as a (loopless) composition $a_1a_2\dots a_M$ of single arcs 
$\{ a_, a_2, \dots, a_M\}$. 
Furthermore, given another loopless presentation 
$D = b_1\dots b_N$ by arcs, these are related by 
repeating one of the fundamental planar identities 
(Figure \ref{isotopy}) locally. 

\begin{figure}[h]
\input{isotopy.tpc}
\caption{\label{isotopy}}
\end{figure}

In the process of applying these identities, 
the composed morphism $F(a_1)\dots F(a_M)$ remains unchanged 
because of the validity of the corresponding 
rigidity identities in $\cT$. 

Thus, for a diagram $D \in K_{m,n}$, 
the morphism $F(D) \in \Hom(Y^{\otimes m},Y^{\otimes n})$ is 
well-defined by the formula 
\[
F(D) = F(a_1)F(a_2)\dots F(a_M),
\]
where $D = a_1a_2\dots a_M$ is a loopless presentation of $D$ 
as a product of single arcs $a_1,a_2,\dots,a_M$. 
By linearity, $F$ is extended to linear maps 
\[
\Hom(X^{\otimes m},X^{\otimes n}) \to 
\Hom(Y^{\otimes m},Y^{\otimes n}).
\]

By the choice of $d$, these linear maps preserve multiplications and 
therefore it gives rise to a functor $\cK_d \to \cT$, 
which, by the construction, is monoidal and satisfies 
$F(\epsilon) = \epsilon_Y$, $F(\delta) = \delta_Y$. 
\end{proof}

\begin{Proposition}
The Temperley-Lieb category $\cK_d$ is a C*-tensor category 
if and only if $d = \pm (q + q^{-1})$ with $q > 0$ or 
$q = e^{i\pi/n}$ ($n \geq 3$). 
For such a value of $d$, the C*-tensor category structure 
on $\cK_d$ is unique up to C*-monoidal equivalences.

An explicit choice is given by 
\[
D^* = \left( \frac{d}{|d|} \right)^{\sharp(D)} D',
\]
where $D'$ denotes the diagram obtained from $D$ up-side down 
(i.e., the reflection of $D$ with respect to a horizontal line)
and $\sharp(D)$, called the arc index of $D$,  
is the difference of the number of $\epsilon$'s and 
$\delta$'s inside $D$. 

\end{Proposition}

\begin{proof}
Assume that $\cK_d$ is a C*-tensor category. Since the star operation 
preserves the central decomposition in a C*-algebra, 
the decomposition $\End(X\otimes X) = \C(1-e) + e$ 
with $e = d^{-1}h$ implies $e^* = e$. Consequently, we have 
$e_j^* = e_j$ for any $j \geq 1$ and the relation 
$d^2 e_1e_2e_1 = e_1$ compels the positivity of $d^2$, i.e.,  
$d \in \R^\times$. 

The inductive analysis on semisimplicity now shows that 
the Jones-Wenzl idempotents $f_k$ are (orthogonal) projections 
and $[2][n]/[n+1] \geq 0$ if $[k] \not=0$ for $k \leq n+1$ 
because 
$[2][n]/[n+1] (f_n\otimes 1_X) e_n (f_n\otimes 1_X)$ is 
a subprojection of $f_{n+1}$. 
The condition is then equivalent to (i) $(d/|d|)^{k+1} [k] > 0$ 
for all $k$ or (ii) $(d/|d|)^{k+1} [k] > 0$ for $1 \leq k <l$ with 
$[l] = 0$ for a positive integer $l \geq 3$. 
It is then immediate to see that these conditions are equivalent to 
the ones given in the statement of the proposition. 

Since both of $\Hom(X\otimes X,I)$ and $\Hom(I,X\otimes X)$ 
are one-dimensional, we should have 
$\epsilon^* = c \delta$ with $c \in \C^\times$ and the positivity of 
$\epsilon^*\epsilon = cd e$ shows that $c$ is a real number, 
which is positive or negative according to the signature of $d$. 
In particular, $\delta^* = c^{-1} \epsilon$. 

Conversely, assume that $d \in \R$ is in the range 
specified above and let $r$ be a positive real. 


For each pair $(m,n)$, define the conjugate-linear map 
$\Hom(X^m,X^n) \to \Hom(X^n,X^m)$ so that 
\[
D^* = \left( \frac{rd}{|d|} \right )^{\sharp(D)} D'
\]
for $D \in K_{m,n}$. 
(We have particularly $\epsilon^* = rd|d|^{-1}\delta$ and 
$h_j^* = h_j$.) 
From the definition, we have $(f\otimes g)^* = f^*\otimes g^*$ 
for morphisms $f$, $g$ in $\cK_d$. 
The map is involutive because of $\sharp(D') = -\sharp(D)$ and 
it satisfies $(fg)^* = g^*f^*$, which follows from 
the additivity of arc index: 
$\sharp(CD) = \sharp(C) + \sharp(D)$ for a composable pair 
$(C,D)$ of diagrams.

In this way, we have defined a *-operation in the tensor category 
$\cK_d$. Note here that only the reality of $r$ and $d$ 
is used up to now. 
Notice also that the Jones-Wenzl idempotents $f_n$ are 
hermitian, i.e., $f_n^* = f_n$, which is checked by the recursive 
formula for them by using the reality of $d$ and hence of $[n]$. 

We shall now check the positivity of the *-operation in question, 
which will be achieved by seeing that there are plenty of 
positive (unitary) representations of $\End(X^n)$. 
By our choice of signature in the definition of 
star operation, we see that the morphism 
\[
r^{1/2} 
\left| 
\frac{\trace(f_{n-1})}{\trace(f_n)} 
\right|^{1/2} 
\varphi_n: X_{n-1} \to X_n\otimes X
\]
gives a realization of $X_{n-1}$ as an orthogonal component 
of $X_n\otimes X$: 
\[
r
\left| 
\frac{\trace(f_{n-1})}{\trace(f_n)} 
\right|
\varphi_n^*\varphi_n 
= \frac{\trace(f_{n-1})}{\trace(f_n)} 
\psi_n\varphi_n = f_{n-1}, 
\]
where, in the first equality, we have calculated as 
\begin{align*}
\varphi_n^* 
&= ((f_n\otimes 1_X)(1_{X^{n-1}}\otimes \delta)^*
= (1_{X^{n-1}}\otimes \delta^*)(f_n\otimes 1_X) \\
&= \frac{|d|}{rd} (1_{X^{n-1}}\otimes \epsilon) (f_n\otimes 1_X)
= \frac{|d|}{rd} \psi_n.
\end{align*}
Notice also that,
by the assumption on $d$, the signature of $\trace(f_k) = [k+1]$ 
alternates as $k$ increases until it vanishes.  


Thus, the path basis in the representation space 
$\Hom(X_k,X^{\otimes n})$ of $\End(X^{\otimes n})$ 
($0 \leq k \leq n$ with $k \equiv n \mod 2$) is orthonormal 
with respect to the inner product $(\cdot|\cdot)$ defined by 
$f^*g = (f|g) f_k$. 

By the obvious decomposition 
\[
\Hom(X^m,X^n) \cong 
\bigoplus_k \Hom(X_k,X^n)\otimes \Hom(X^m,X_k),
\]
the previous observation on representations of 
$\End(X^{\otimes n})$ reveals that the category $\cK_d$ turns out to 
be a C*-tensor category. 

Finally, we show the uniqueness of C*-structures. 
Let $*$ denote a star operation in $\cK_d$ which makes $\cK_d$ 
into a C*-tensor category. 
Given $\lambda \in \C^\times$, define a monoidal functor 
$F: \cK_d \to \cK_d$ so that 
\[
F(\epsilon) = \lambda \epsilon, 
\quad 
F(\delta) = \lambda^{-1} \delta,
\]
i.e., 
\[
F(D) = \lambda^{\sharp(D)} D
\]
for a diagram $D$. 
Then $F$ is an automorphism of the tensor category $\cK_d$ 
and the new star operation $\star$ in $\cK_d$, 
which again gives a compatible C*-structure in $\cK_d$, is defined 
by the relation 
\[
F(f)^\star = F(f^*)
\quad
\text{for a morphism $f$ in $\cK_d$.}
\]
Then, on the level of diagrams, we have the explicit formula 
\[
D^\star = |\lambda|^{-2\sharp(D)} D^*.
\]

Thus, compatible structures of C*-tensor category on $\cK_d$ 
are unique up to C*-monoidal equivalences.
\end{proof}

\begin{Remark}
The duality isomorphism $X \cong X^{**}$ is given by the identity 
morphism $1_X$ for the case $d>0$, whereas it is $-1_X$ if $d<0$, 
resulting in the positive value $|d|$ as 
the quantum dimension of $X$ in either cases. 
A bit more analysis (see \cite{FK}, \cite{KW} for example) shows that
$\cK_d$ and $\cK_{-d}$ are related to each other 
by twisting associativity isomorphisms in monoidal structure 
with respect to 
a non-trivial 3-cocyle of the group $\Z^2$. 
On the level of operator algebras, 
this is interpreted as twisting a generating bimodule 
with respect to an outer automorphism $\alpha$ of a factor $N$ 
satisfying (i) $\alpha\circ \alpha = \text{Ad}\,u$ (ii) and 
$\alpha(u) = -u$ with $u$ a unitary in $N$. 
\end{Remark}

\appendix
\section{Universality Property}

We shall here present a proof of the fact that 
the Kauffman's planar algebra $\C[K_n,d]$ is 
identified with the Temperley-Lieb algebra $A_n$, 
which is universally generated by elements 
$\{ h_1, \dots, h_{n-1}\}$ with the relations of Temperley-Lieb.
Since there is a natural homomorphism 
$\pi: A_n \to \C[K_n,d]$ by universality, 
the problem is in checking the bijectivity of $\pi$, 
which was claimed in \cite{Kau} 
with more accounts supplied 
in \cite[Theorem 4.3]{Kau2} (cf.~\cite{BD} also).
The proof obviously consists of two parts: 
the surjectivity of $\pi$ and the injectivity of $\pi$. 
The former is the generating property of elementary diagrams 
in the algebra $\C[K_n,d]$, while the latter is reduced to 
the problem of counting reduced words of generators.

As for the generating property, 
the following would not be the shortest proof, 
compared with the one 
given in \cite[Theorem XII.3.2]{Tu} for instance,  
but has the advantage 
that it produces reduced words. 
(The Jones' normal form is then obtained shearing positions 
of mutually commutable elementary diagrams as indicated by 
\cite[Figure 16]{Kau2}. 
  
For discussions of the proof, we introduce some terminologies first. 
Given a diagram $D$ in $K_n$, a string inside $D$ 
is called a through string if it connects 
upper and lower vertices, and an arc otherwise.
A through string is said to be vertical if it connects 
vertices in the same horizontal position.
A handle is, by definition, 
an arc which connects neighboring vertices.   

We shall apply an induction of trying to increase the number of handles 
inside relevant diagrams to get the generating property.
If a diagram $D$ contains a vertical through string, then 
$D$ is of the form $D'\otimes D''$ with $D' \in K_{n'}$ and 
$D'' \in K_{n''}$, whence the problem is reduced 
to diagrams of less strings.

Consider the case 
that $D$ contains a through string connecting two vertices 
$i < j$, say $i$ on the top and $j$ on the bottom. 
Then $j-i$ is an even number (otherwise there would appear 
unconnected vertices) and we can apply waving to the string 
(i.e., the string is deformed so that it repeats local maxima and minima 
alternately)
after separating it from the other strings by stretching out 
sufficiently. 
Then there arise three patterns 
depending on the position where the vertex $i+1$ on the top 
is terminated: 
If it ends at a bottom vertex numbered by $k$, then $k > j$ 
and the waving of this second through string, 
together with the waving of the first string, gives us 
couplings of handles during the horizontal interval $[i, j]$. 
If we further stretch the waving of the second string 
on the unoverlapping interval $[j+1,k]$, 
push down the strings tied to vertices in the interval 
$[1,i-1]$, and pull up the strings ending at vertices in the interval 
$[k+1,n]$ sufficiently enough 
so that $D$ becomes the composition of three diagrams in $K_n$ 
as indicated by dotted horizontal lines in Figure \ref{waving1}. 
Then the middle diagram is apparently a product 
of (mutually commuting) elementary diagrams, whereas 
the remaining diagrams contain vertical through strings. 

Otherwise, the string starting from $i+1$ forms an arc ending 
at the vertex $k$ with $k > i+1$. Though we need to further 
divide into two patterns depending 
on the relative position of $j$ and $k$, 
a similar decomposition is possible as indicated by 
pictures (Figure \ref{waving2}, Figure \ref{waving3}) 
and we are again reduced to diagrams of less strings. 

Finally, there remains the case that $D$ contains no through 
string. If there is an arc which is not a handle, 
we see $D$ containing a part of continuing handles 
surrounded by one arc. 
Then, by waving the surrounding arc, and a similar rearrangement
as above allows us to identity $D$ with 
the composition of two diagrams such that 
one of them is again a product of (commuting) elementary diagrams 
and the other has more handles than the original diagram 
(Figure \ref{waving4}). 
Repeating the same procedure, we end up with diagrams in which 
all the arcs are handles, thus again a product of commutating 
elementary diagrams.  

\begin{figure}[h]
\input{wavingo1.tpc}  
\caption{\label{waving1}}
\end{figure}

\begin{figure}[h]
\input{wavingo2.tpc}  
\caption{\label{waving2}}
\end{figure}

\begin{figure}[h]
\input{wavingo3.tpc}  
\caption{\label{waving3}}
\end{figure}

\begin{figure}[h]
\input{wavingo4.tpc}  
\caption{\label{waving4}}
\end{figure}

To see the injectivity of the map $\pi: A_n \to \C[K_n]$, 
we use the dimension estimate of $A_n$ here:
According to V.~Jones, 
by a {\bf word} in $A_n$, we shall mean a product 
$h_{i_1}h_{i_2}\dots h_{i_k}$ where $h_{i_1},\dots, h_{i_k} 
\in \{ h_1,\dots, h_{n-1} \}$ with 
two words identified if we can relate them each other  
by applying the commutativity $h_ih_j = h_jh_i$ ($|i-j| \geq 2$) 
to their ingredients. Thus $h_1h_3h_1h_2 = h_3h_1h_1h_2$ 
as a word for example. 
The length of a word is, by definition, the number of 
$h_i$'s appearing in the word. 
A word is said to be {\bf reduced} if its length is minimal under 
the replacements of $h_ih_{i\pm 1}h_i$ to $h_i$ and $h_i^2$ to $h_i$.
Thus any word is equal to a reduced one up to multiplication of 
powers of $d$. 

From the commutation relations, $h_m$ with $m$ the maximal index 
appears only once in a reduced word. 
According to V.~Jones, a reduced word is further relocated so that 
$h_m$ is placed at the rightest end. 
Then, after the point $h_m$, there follows  
a sequence of the form $h_m h_{m-1} \dots h_l$ with $l \leq m$. 
On the left of the block of this sequence, we are left a reduced word 
consisting of elements in $\{ h_1, h_2, \dots, h_{m-1}\}$, 
for which we can apply the same procedure to get eventually the form 
\[
(h_{i_1}h_{i_1-1}\dots h_{j_1}) 
(h_{i_2}h_{i_2-1}\dots h_{j_2}) 
\dots
(h_{i_k}h_{i_k-1}\dots h_{j_k}), 
\]
with $i_1 \geq j_1, i_2 \geq j_2, \dots, i_k \geq j_k$ and 
$1 \leq i_1<i_2<\dots<i_k \leq n-1$
(the case $k=0$ corresponding to the empty word). 

Since this is assumed to be a reduced word, we should have 
$1 \leq j_1<j_2<\dots<j_k \leq n-1$ as well. 
In fact, if $j_1 \geq j_2$ for example, we should have 
$h_{j_1}$ appearing in the block 
$h_{i_2}h_{i_2-1}\dots h_{j_2}$ by $j_1 \leq i_1 < i_2$ and hence 
a reduction of the form $h_{j_1}h_{j_ 1 + 1} h_{j_1} = h_{j_1}$ 
takes place, reducing the word length.

\begin{Proposition}[{\cite[\S 4]{Jo}}]
Any word in $A_n$ is equal to a reduced one up to scalar 
multiplications and
the number of reduced words in $A_n$ is given by the Catalan number 
\[
C_n = 
\frac{1}{n+1} 
\begin{pmatrix} 2n\\ n\end{pmatrix}. 
\]
\end{Proposition}

\begin{Corollary}[cf.~{\cite[Theorem 4.3]{Kau2}}]
The algebra $\C[K_n,d]$ is 
universally generated by 
$\{ h_1, \dots, h_{n-1} \}$ and the unit $1$ with the relations 
\[
h_i^2 = d h_i, 
\quad 
h_ih_j = h_jh_i 
\quad (|i-j| \geq 2), 
\quad 
h_i h_{i\pm 1}h_i = h_i
\]
for $1 \leq i,j \leq n-1$, whence it is identified with 
the Temperley-Lieb algebra $A_n$. 

In particular, for $n \geq 1$, 
the obvious homomorphism $A_n \to A_{n+1}$ of 
Temperley-Lieb algebras is injective. 
\end{Corollary}

In connection with the identification $A_n = \C[K_n,d]$, 
the following answers how Jones' reduced words are 
related to Kauffman's diagrams, which also constitutes a part of 
\cite[Theorem 4.3]{Kau2}. We shall present a proof as a continuation 
of discussions so far. 

\begin{Proposition}
The set $K_n$ of Kauffman diagrams is exactly 
the image of the set of reduced words in $A_n$ 
under the natural isomorphism.
\end{Proposition}

\begin{proof}
We need to show that reduced words contain no loops when 
they are computed as compositions of diagrams in $K_n$. 
We shall check this by an induction on the number of blocks 
in the Jones normal form. 
Let $h_mh_{m-1}\dots h_l$ be the last block in a Jones reduced word. 
As a diagram, this descending sequence is given by 
Figure \ref{descend}. 
Since the previous blocks constitute a reduced word of Jones form 
with the number of blocks reduced, it contains no loops inside 
by the induction hypothesis and therefore the composition with 
the block $h_m \dots h_l$ remains loopless. 
\end{proof}

\begin{figure}[h]
\input{descendo.tpc}
\caption{\label{descend}}
\end{figure}



\begin{thebibliography}{8}
\bibitem{Ba}
T.~Banica, Quantum groups and Fuss-Catalan algebras,
arXiv:math.QA/0010084.
\bibitem{BW} 
J.W.~Barrett and B.W.~Westbury, 
Spherical categories, 
{\it Adv.~in Math.}, 
143(1999), 357--375. 
\bibitem{BJ}
D.~Bisch and V.~Jones, 
Algebras associated to intermediate subfactors, 
{\it Invent.~math.}, 
128(1997), 89--157.
\bibitem{BD}
M.~Borisavljevi\'c and K.~Dosen, 
Kauffman monoids, 
\textit{J.~Knot Theory and Its Ramifications}, 
11(2002), 127--143.
\bibitem{CP}
V.~Chari and A.~Pressley, 
{\it A Guide to Quantum Groups}, 
Cambridge Univ.~Press, 
1995.
\bibitem{Da}
A.A.~Davydov, 
Monoidal categories,  
{\it J.~Math.~Sci.~(New York)}, 
88(1998), 457--519.
\bibitem{EK}
D.~Evans and Y.~Kawahigashi, 
\textit{Quantum Symmetries on Operator Algebras}, Clarendon Press, 
Oxford, 1998. 
\bibitem{FK}
J.~Fr\"ohlich and T.~Kerler, 
\textit{Quantum Groups, Quantum Categories and Quantum Field Theory}, 
Lec.~Notes in Math.~1542, Springer-Verlag, 1993. 
\bibitem{GHJ}
F.M.~Goodman, P.~de la Harpe and V.F.R.~Jones, 
{\it Coxeter Graphs and Towers of Algebras}, 
Springer-Verlag, 1989. 
\bibitem{GW}
F.M.~Goodman and H.~Wenzl, 
Ideals in the Temperley-Lieb category, 
arXiv:math.QA/0206301.
\bibitem{Jo}
V.F.R.~Jones, Index for subfactors, 
{\it Invent.~Math.}, 
72(1983), 1--25.
\bibitem{Kas}
C.~Kassel, 
{\it Quantum Groups}, 
Springer-Verlag, 
Berlin-New York, 
1995. 
\bibitem{Kau}
L.H.~Kauffman, State models and the Jones polynomial, 
{\it Topology}, 
26(1987), 395--407.
\bibitem{Kau2}
\underline{\phantom{Kauffman}}, 
An invariant of regular isotopy, 
\textit{Trans.~Amer.~Math.~Soc.}, 
318(1990), 417--471.
\bibitem{Kau3}
\underline{\phantom{Kauffman}}, 
Knots and Physics, 1991, World Scientific Publishing.
\bibitem{KW}
D.~Kazhdan and H.~Wenzl,  
Reconstructing monoidal categories, 
{\it Adv.~Soviet Math.}, 
16(1993), 111--136.  
\bibitem{GKP}
R.L.~Graham, D.E.~Knuth and O.~Patashnik, 
\textit{Concrete Mathematics}, 
Addison-Wesley, 1994.
\bibitem{La}
Z.A.~Landau, Fuss-Catalan algebras and chains of 
intermediate subfactors, 
{\it Pacific J.~Math.}, 
197(2001), 325--367.
\bibitem{Mac}
S.~MacLane, 
{\it Categories for the Working Mathematician}, 
Springer-Verlag, 
Berlin-New York, 
1971. 
\bibitem{Mal}
G.~Maltsiniotis, 
Traces dans les cat\'egories monoidales, dualit\'e et cat\'egories monoidales fibries, 
{\it Cahiers Topologie G\'eom. Diff\'erentielle Cat\'eg.},
36(1995), 195--288.
\bibitem{PP}
M.~Pimsner and S.~Popa, 
Iterating the basic constructions, 
{\it Trans.~Amer.~Math.~Soc.}, 310(1988), 127--134.
\bibitem{TL}
H.~Temperley and E.~Lieb, Relations between the 
`percolation' and `coloring' problem and 
other graph-theoretical problems associated with 
regular plane lattices, 
{\it Proc.~Roy.~Soc.~London A}, 322(1971), 251--280.
\bibitem{Tu}
V.G.~Turaev, 
\textit{Quantum Invariants of Knots and 3-Manifolds}, 
Walter de Gruyter, Berlin-New York, 1994. 
\bibitem{W}
B.W.~Westbury, 
The representation theory of the Temperley-Lieb algebras, 
\textit{Math.~Z.}, 219(1995), 539--565. 
\bibitem{NOA}
S.~Yamagami, 
A note on Ocneanu's approach to Jones' index theory, 
{\it Internat.~J.~Math.}, 
4(1993), 859--871. 
\bibitem{FTC}
\underline{\phantom{S.~Yamagami}}, 
Free products of semisimple tensor categories, preprint.
\end{thebibliography}
\end{document}